\documentclass[a4paper,10pt]{article}
\def \Z {{\mathbf {Z}}}

\def \N {{\mathbf {N}}}

\def \B {\mathcal {B}}

\def\eps{\varepsilon}
\title{  Cлабо   замкнутые  дискретные полугруппы, \\
       порожденные  фукциями  от унитарного оператора}
\author{В.В.Рыжиков} 
\date{}
\usepackage[T1]{fontenc} % иногда T1 вместо T2A
\usepackage[cp1251]{inputenc}  
\usepackage[russian]{babel}
\usepackage{graphicx}
\graphicspath{{pictures/}}
\DeclareGraphicsExtensions{.pdf,.png,.jpg}
\begin{document}
\Large
\maketitle
\begin{abstract}  
\bf Weakly closed discrete semigroups generated by given operator-valued functions. \rm
The function $P(T)=\sum_{i=0}^\infty c_i T^i$  is admissible if
$c_i\geq 0$, $\sum_{i=0}^\infty c_i\leq 1$.
For any given set of admissible functions $P_1,\dots, P_k$  there is  a unitary operator $T$ of dynamic origin such that the weak closure of its powers  is a semigroup generated by the operators $0$, $T$, $P_1(T),\dots, P_k(T),T^\ast, P_1(T^\ast ),\dots, P_k(T^\ast)$.

Функция $P(T)=\sum_{i=0}^\infty c_i T^i$, от  унитарного оператора $T$  допустима, если 
$c_i\geq 0$, $\sum_{i=0}^\infty c_i\leq 1$.
Для заданного набора допустимых функций  $P_1,\dots, P_k$  доказано существование унитарного оператора $T$ динамического происхождения, который имеет непрерывный спектр, а слабое замыкание индуцированного им действия   является полугруппой, порожденной операторами $0$, $T$,   $P_1(T),\dots, P_k(T)$ и их сопряженными.
\rm

%Библиография: 10 названий, УДК: 517.987, \ MSC: Primary 28Y05; Secondary 58F11

Ключевые слова и фразы: \it  эргодическое действие,
слабое замыкание, функция от оператора,  полугруппа операторов.\rm

\end{abstract} 

\section{Введение} 
Слабые замыкания бесконечномерных групп в унитарных представлениях являются источником  разнообразных полугрупповых  структур   \cite{Ol}. Эти же структуры возникают в замыкании действий на пространствах с мерой \cite {N}. Cлабые замыкания сохраняющих меру действий находят ряд  применений, из недавних работ см.  \cite{R19}--\cite{R22}. Настоящая работа посвящена задаче реализации слабого замыкания  эргодического действия в виде полугруппы с заданной алгебраической  структурой. 
 
Простейшей структурой слабого замыкания  обладают так называемые  перемешивающие действия.  Для них слабое  замыкание  состоит в добавлении  к действию  одной предельной  точки.
%В случае действия на  вероятностном  пространстве этой точкой является % $\Theta$ --  оператор ортопроекции на пространство констант. Для %действий на пространстве с сигма-конечной мерой таким является нулевой %оператор. 
 Доказательство   свойства перемешивания   может оказаться нетривиальной  задачей (см., например, \cite{O},\cite{SH}).
Описать поностью   слабое  замыкание  для неперемешивающих действий  в общем случае  еще сложнее.  Даже  поиск некоторых нетривиальных слабых пределов  может потребовать привлечения тонких методов (см. \cite{K1}).

  Наша цель   -- предъявить      эргодическое $\Z$-действие, для которого слабое замыкание является  полугруппой,  порожденной заданными  функциями $P(T)$ от  унитарного оператора $T$.  Пусть для ряда  $P(T)=\sum_{i=0}^\infty c_i T^i,$   выполнены условия:
   $\sum_{i=0}^\infty c_i  \leq 1,$ $c_i\geq 0.$  
Для таких функций  $P$ значения    $P(T)$  могут быть слабыми пределами степеней унитарного оператора $T$.  Мы дополнительно потребуем, чтобы   $c_0>0$ и  для некоторого $k>0$  выполнялось $c_k>0$. Это условие обеспечит слабую сходимость 
$P(T)^n\to_w 0$, $n\to\infty$. Спектр  такого оператора $T$   непрерывен.  Функции $P(T)=\sum_{i=0}^\infty c_i T^i,$ с указанными ограничениями на коэффициенты $c_i$  назовем  \it допустимыми.  \rm 

Для всякого набора допустимых функций  $P_1,\dots, P_k$  будет доказано существование автоморфизма $T$ пространства с сигма-конечной мерой, для которого слабое замыкание его степеней    является полугруппой, порожденной операторами $0$, $T$,  $P_1(T),\dots, P_k(T)$ и операторами, сопряженными к ним.

Слабое замыкание унитарного действия  $\{T^i:\,i\in \Z\}$ содержит все степени оператора $T$  и все предельные точки степеней, т.е.   операторы $P$, для которых  выполнено 
$T^{i_k}\to_w P$ при ${i_k}\to\infty$. Множество  таких пределов  обозначаем через $Lim(T)$. В качестве примера отметим, что в \cite{R19}  описан класс   обратимых преобразований $T$  с сигма-конечной инвариантной мерой, для которых 
 $$Lim(T)=\{0\}\cup \{2^{-n}{T^z}\, :  \, n\in\N,\, z\in\Z \}.$$ 
Для  автоморфизма  $T$ вероятностного   пространства   свойство перемешивания 
означает слабую сходимость $$T^i\to_w \Theta, \ i\to\infty, $$
где $\Theta$ -- оператор ортопроекции $L_2$ на пространство констант.
Преобразование и  отвечающий ему унитарный оператор в статье обозначаются одинаково, они связаны соотношением   $Tf(x)=f(T(x))$. Если ограничить оператор $T$ на пространство, ортогональное константам, то свойство перемешивания выражается сходимостью 
$$T^i\to_w 0, \ i\to\infty.$$
Сходимость $T^i\to_w 0$, $i\to\infty$,   является определением перемешивания  для автоморфизма $T$ пространства с сигма-конечной мерой. Темин "перемешивание" в  этом случае, вообще говоря, не имеет отношения к физическому перемешиванию. Обычный сдвиг на прямой  является перемешивающим. Отчасти терминологию оправдыват пуассоновские (и гауссовские)  надстройки $T_\circ$ (см. \cite{R14},\cite{PR}), превращая номинальное перемешивание $T^i\to_w 0$ в настоящее: $T_\circ^i\to_w\Theta$.

Обозначим через  $\langle T, Q_1,$ $\dots,$ $Q_k\rangle$ полугруппу, порожденную операторами $0$, $T, Q_1$,$\dots$, $Q_k$, и операторами,  сопряженными к ним. 
Основной результат статьи формулируется следующим образом.

\vspace{3mm}
\bf  Теорема 1. \it Для всякого набора допустимых функций  $P_1,\dots, P_k$  существует унитарный оператор $T$ с непрерывным  спектром такой, что   
$$Lim(T)=\langle P_1(T),\dots, P_k(T)\rangle\setminus\{T^n: n\in\Z\}. $$\rm

\vspace{3mm}
Искомый   оператор $T$  индуцирован преобразованием, сохраняющим сигма-конечную меру. Он ищется в виде конструкции, которая   задается счетным набором параметров. Часть   параметров  обеспечивает нужные слабые пределы степеней преобразования. Другая часть параметров играет роль фильтра, избавляя слабое замыкание от  ненужных   пределов.

Этот результат имеет приложение к теории автоморфизмов вероятностного пространства. Гауссовские и пуассоновские надстройки над $T$ унитарно эквивалентны оператору
$e^\otimes T:=\bigoplus_{n=0}^\infty  T^{\odot n}$, поэтому    
 $$Lim(e^\otimes T)=\langle e^\otimes P_1(T),\dots, e^\otimes P_k(T)\rangle.$$ 
Из  свойств гауссовских и пуассоновских надстроек вытекает, что в нашем случае они неизоморфны как автоморфизмы: централизатор первых континуален, а централизатор вторых счетен (см. \cite{R19},\cite{PR}  для обоснования этого контраста).

\section {Конструкции преобразований}
 Нам понадобятся   конструкции  сохраняющих меру преобразований  ранга один. Напомним их определение.
 Пусть задано натуральное число $h_1$  и  последовательность натуральных чисел $r_j\to\infty$ 
вместе с  последовательностью  целочисленных векторов (параметров надстроек) 
$$ \bar s_j=(s_j(1), s_j(2),\dots, s_j(r_j-1),s_j(r_j)), \ s_j(i)\geq 0.$$ 
Параметры $h_1$, $r_j$, $\bar s_j$ полностью  определяют конструкцию преобразования.   Опишем ее  построение. 

На   шаге $j\geq 1$ определена 
система  непересекающихся полуинтервалов, одинаковой длины (башня высоты $h_j$)
$E_j, TE_j, T^2E_j,\dots, T^{ h_j-1}E_j.$
На этих полуинтервалах, кроме последнего,
пребразование $T$ действует как обычный  перенос полуинтервалов. 
На  $T^{ h_j-1}E_j$ оно пока не определено. 

 Представим 
 $E_j$ в виде  дизъюнктного объединения  полуинтервалов $ E_j^i$, $1\leq i \leq r_j$, одинаковой меры. Набор  
$E_j^i, TE_j^i ,T^2 E_j^i,\dots, T^{h_j-1}E_j^i$ называется 
$i$-ой колонной на этапе $j$.
К этому набору добавим $s_j(i)$ полуинтервалов меры $\mu(E_j^i)$, 
 тем самым получив надстроенную колонну: набор непересекающихся полуинтервалов
$$E_j^i, TE_j^i, T^2 E_j^i,\dots, T^{h_j-1}E_j^i, T^{h_j}E_j^i, T^{h_j+1}E_j^i, \dots, T^{h_j+s_j(i)-1}E_j^i.$$
Терерь соберем   надстроенные колонны в  башню. Для этого при $i<r_j$ положим $T^{h_j+s_j(i)}E_j^i = E_j^{i+1}.$ 
Обозначим  $E_{j+1}= E^1_j$. Под действием степеней $T^n$ при $n<h_{j+1}$, где  
$$ h_{j+1} =h_jr_j +\sum_{i=1}^{r_j}s_j(i),$$
полуинтевал $E_{j+1}$ пробежит этажи всех надстроенных колонн. Таким образом, мы получили  башню этапа $j+1$:
$$E_{j+1}, TE_{j+1}, T^2 E_{j+1},\dots, T^{h_{j+1}-1}E_{j+1}.$$

Продолжая построение до бесконечности,  мы определяем  преобразование $T$   на объединении  $X$ всех рассматриваемых полуинтервалов. 
Оно обратимо и  сохраняет обычную меру Лебега на интервалах. 
Известно, что $T$  эргодическое и имеет простой спектр. 
Приведем примеры конструкций и некоторые  их свойства. 

\vspace{2mm}
\bf  \ 1. \bf Параметры конструкции $T$: $\bf r_j\to\infty$, $\bf s_j(i)=h_j$. \rm Для соответствующего  оператора  $T$ наблюдаем любопытный эффект:
 $$  T^{h_j} \to_w 0, \ \   T^{2h_j} \to_s I.$$ 
Действительно, если множества $A,B$ состоят из этажей башни на этапе
$j_0$, то для всех $j\geq j_0$  имеем
$$ T^{h_j}A\cap B= \emptyset,   \ \  \mu (A\Delta T^{2h_j}A)= 2\mu(A)/r_j.$$ Эти равенства очевидным образом обеспечивают нужные сходимости. Для гауссовских и пуассоновских надстроек $S$ над $T$ мы имеем 
$$  S^{h_j} \to_w \Theta, \ \   S^{2h_j} \to_s I.$$

\vspace{2mm}
\bf  2. \ $\bf  Lim(T)=\{0\}\cup \{2^{-n}{T^z}\}.$ \rm
Параметры конструкции $T$: $r_j=2$, $s_j(1)=0$, $s_j(2)=j h_j$.  Такой автоморфизм $T$ обладает тривиальным централизатором. См. \cite{R19}.

\vspace{2mm}
\bf  3. Все допустимые пределы. \rm  Параметры автоморфизма $T$, для которого $ Lim(T)$ содержит  всевозможные $ P(T)$  для всех допустимых функций $P$: $$r_j=3, \ s_j(1)=h_j, \ s_j(2)= h_j+ [\sqrt j],\ s_j(3)=h_j.$$
Аналогичный пример рассмотрен в \cite{R20}, теорема 5.1.

\vspace{2mm}
\bf  4. Стохастический ансамбль. \rm В статье \cite{O} не предъявлялась конкретная перемешивающая конструкция ранга один,  а вероятностными методами доказывалось существование класса таких  конструкций.
Для нахождения  автоморфизма $T$ с полугруппой 
$$Lim(T)=\langle P\rangle, \ P={(I+T)}/{2}, \eqno (1)$$ вероятностный подход представляется весьма естественным, так как для случайной последовательности $s_j(i)$  для фиксированных $a_k\in \{0,1\}$ частота события   " $s_j(i+k)=a_k$ при $k=1,2,\dots, m$ "
 стремится к $2^{-m}$. А это обеспечивает сходимость
$T^{-mh_j}\to_w P^m$, что  нам и нужно.   

Положим  $r_j=2^j$,  $s_j(i)\in \{h_j,h_j+1\}$. Для фазового пространства $X$ такой конструкции имеем  $\mu(X)=\infty$.
На множестве всех таких конструкций очевидным образом вводится равномерная мера ($s_j(i)$ принимает независимо значения $h_j$ и $h_j+1$ с вероятностью $1/2$).

Гипотеза: \it  для почти всех  указанных конструкций   выполняется $(1)$.\rm
\\
Основная задача -- обеспечить отсутствие  ненужных слабых пределов.  Мы  не проводим   статистический  анализ, а  модифицируем конструкцию, гарантированно избавляясь от лишних пределов. С этой целью  будет рассмотрим класс   конструкций,  для которых при  $i$,  не кратном $j$, определены  небольшие параметры $s_j(i)$ случайной природы. Для  $i$, кратных $j$, , большие значения   $s_j(i)$   задаются, например так:
$$jh_j<s_j(j), \ \  js_j(kj)<s_j(kj+j).$$
Небольшие параметры $s_j(i)$ обеспечат наличие в $Lim(T)$  полугруппы $\langle P_1(T),$ $\dots,$ $ P_k(T)\rangle$, а большие параметры позволят установить, что слабое замыкание изучаемого действия не содержит  других предельных точек. 

Идея использовать модификации  больших параметров
была навеяна так называемыми сидоновскими конструкциями.  В случае пространства с бесконечной мерой
благодаря полной свободе в выборе параметров  надстроек  несложно получить перемешивающую конструкцию ранга один. В \cite{R14} были рассмотрены так называемые сидоновские надстройки, для которых свойство перемешивания  вытекало из определений. Например, при выборе параметров надстроек, удовлетворяющих условию
$$h_j\ll s_j(1)\ll s_j(2)\ll\dots\ll  s_j(r_j),\ r_j\to\infty,$$
для всяких множеств $A,B\in\B$ конечной меры для всех больших $m\in[h_j,h_{j+1})$ выполняется
$$\mu (T^mA\cap B) \leq {(\mu(A)+\eps)} /{r_j}.$$

\section{ $P$-конструкции} 
Пока рассмотрим  задачу построения преобразования  $T$ такого, что  
$$Lim(T)=\langle P(T)\rangle,  \ \ P(T)=\sum_0^n c_kT^k, \ \ \sum_0^n c_k=1. $$  (На самом деле общий случай не имеет принципиальных отличий
от этого частного.)    Пусть $r_j\to\infty$,  выбираем параметры 
$s_j(i)\in\{0,1,\dots, n\}$  случайно с вероятностью события  $s_j(i)=k$, равной $c_k$. Тогда  для   типичной последовательности  $s_j(i)$, для всякого $m>0$ будет выполнено 
 $$ T^{-mh_j}\to_w P(T)^m.$$
Указанную слабую сходимость  мы будем также записывать 
в виде 
 $$ T^{-mh_j}\approx_w P(T)^m$$
или 
 $$ dist_w(T^{-mh_j}\, ,\, P(T)^m)\to 0.$$
Для $m=1$ имеем
$$T^{-h_j}\approx_w \frac 1{r_j-1}\sum_{i=1}^{r_j-1}T^{s_j(i)}= 
\sum_{k=0}^{\infty}c_{k,j}^{(1)} T^{k}\to_w P(T),$$
так как частота $c_{k,j}^{(1)}$ события $s_j(i)=s$  сходится к вероятности $c_k^{(1)}=c_k$ при $j\to\infty$.  Обозначим 
$P=P(T)$.   C учетом   независимости случайных  значений $s_j(i)$  и $s_j(i+1)$ получаем  
$$T^{-2h_j}\approx_w \frac 1{r_j-2}\sum_{i=1}^{r_j-2}T^{s_j(i)+s_j(i+1)}
=  \sum_{k=0}^{\infty}c_{k,j}^{(2)}T^{k}\to_w P^2(T),$$ 
так как частота $c_{k,j}^{(2)}$ события  $s_j(i)+s_j(i+1)=s $ сходится 
к $c_{k}^{(2)}$ -- коэффициенту при  $T^s$ в ряде $P^2(T)$.
Теперь сформулируем  более общее утверждение.

\vspace{3mm}
\bf Лемма 1. \it Найдутся последовательности $r_j\to\infty$ и  
$\bar s_j$ такие, что  для соответствующей конструкции $T$ и  некоторой последовательности $\eps_j\to 0$  для всех $j$ и  $m=1,2,\dots, j$  выполнено $dist_w(T^{-mh_j},\,  P(T)^m)<\eps_j$.
\rm

\vspace{3mm} 
Доказательство. Существование нужной последовательности параметров $s_j(i)$  устанавливается при помощи эргодической теоремы (иначе говоря, усиленного закона больших чисел), примененной к $S$  --
 схеме Бернулли с образующим разбиением $\{C_0, C_1\dots,  \}$, $\nu(C_k)=c_k$, $\sum_k c_k=1$. 
Рассмотрим функцию  $f$, которая  на  $C_k$ принимает значение $k$. 
Пусть $c^m_k $ -- коэффициент при степени $T^k$ в полиноме $P(T)^m$.  Заметим, что 
$$\nu(Y^m_k) = c^m_k, \ Y^m_k=\{x: \, \sum_{i=1}^{m} f(S^ix)=k\}.$$ 
Для почти всех $x$ в силу эргодической теоремы выполнено 
$$\frac 1{r_j}\sum_{i=1}^{r_j}\chi_{Y}(S^ix)\to \nu(Y).$$
Для $\eps_j$ найдем $r_j$ такое, что для всех $m=1,2,\dots, j$ при  $c^m_k>0$ выполнено
$$\left|c^m_k -\frac 1{r_j}\sum_{i=1}^{r_j}\chi_{Y^m_k}(S^ix)\right|<\eps_j c^m_k.$$
Положим $$s_j(i)= f(S^ix), \ 1\leq i \leq r_j,$$
для некоторой типичной точки $x$.
В силу указанных условий  для всех $m$, $m=1,2,\dots,j$ при $c^m_k>0$ имеем
$$c^m_k -|\{i\,:   1\leq i\leq r_j, \
 s_j(i)+\dots +s_j(i+m-1)=k\}|/r_j \ < \  \eps_j c^m_k.$$

 Таким образом, для   $m=1,2,\dots, j$ будет выполнено
$$T^{-mh_j}\approx_w \frac 1{r_j-m}\sum_{i=1}^{r_j-m}T^{s_j(i)+
\dots +s_j(i+m-1)}= 
\sum_{i=1}^{\infty}c_{i,j}^{(m)}T^{i}\approx_s P^m(T),$$
$$ w=dist(T^{-mh_j},\, P^m(T))<\eps_j.$$

\vspace{3mm} 
\bf $P$-конструкция.  \rm  Параметры конструкции, которую нам предоставляет лемма 1, изменим $s_j(i)$ только для тех  $i$, которые кратны $j$.  Пусть измененные параметры удовлетворяют условиям
$$jh_j<s_j(j), \ \ js_j(kj)<s_j(kj+j). \eqno (Sidon)$$
Полученное преобразование назовем   $P$-конструкцией. \rm

\vspace{3mm}

\bf Теорема 2. \it Если для   $P$-конструкции $T$ оператор $Q\neq 0$
лежит в слабом замыкании степеней $T^m$, то найдутся последовательность ${m_p}$ и натуральное $n$ и целые числа  $a_{1},a_{2},\dots, a_{n}, z$  такие, что для всякого $p$ для некоторых $j_1,j_2,\dots,j_n$
$$m_p= a_{1}h_{j_1}+ a_{2}h_{j_2}+\dots+ a_{n}h_{j_n}+ z, \ j_1>j_2>\dots>j_n,$$
и для некоторых $n_1,n_2$ выполнено
$$T^{m_p}\to_w Q=P^{\ast\,n_1} P^{n_2} T^z.$$ 
\rm

\newpage
%\vspace{3mm}
\bf Лемма 3. \it Пусть $T$ -- унитарный оператор  с непрерывным спектром.  Для всякого набора допустимых функций $P_1,\dots, P_k$ 
при $$Q_i\in \{P_1(T),\dots, P_k(T), P_1(T^\ast),\dots, P_k(T^\ast)\}.$$ 
  имеет место сильная операторная сходимость
$$ \prod_{i=1}^n  Q_i(T)\,\to_s\, 0, \  n\to\infty. $$
В частности, $P(T)^n\to_s 0$ для всякой допустимой функции $P$. 
\rm

\vspace{3mm} 
Доказательство.   Заметим, что $|P_i(z)|<1$ 
для  почти всех $z$ ($|z|=1$) относительно спектральной меры $\sigma$ оператора $T$, так как  равенство $|P_i(z)|=1$ возможно лишь для $z$, являющихся корнями из 1 (напомним, что для ряда $P_i(z)=\sum_kc_kz^k$ выполняется $c_0>0$,  $\sum_k|c_k|\leq 1$).
Поэтому для $\varphi\in L_2(\sigma)$ имеет место сходимость
$$ \left\|\varphi \prod_{i=1}^N  Q_i\right\|\ \to\, 0, \ N\to\infty, $$
что равносильно утверждению леммы.

\vspace{2mm}
\bf Лемма 4. \it Пусть   $T$ является $P$-конструкцией. Для любого $\eps>0$  найдется $L$ такое, что при $1\leq n \leq  L$
выполнено $dist_w(T^{nh_j},  P(T^\ast)^n)<\eps,$ а при $Lh_j<m< h_{j+1}$ имеет место  неравенство $dist_w(T^m,  0)<\eps.$
\rm

\section{ Конструкции с заданной полугруппой слабых пределов} 
Для каждой отдельной функции $P_q$,\, $q=0,1,\dots, k-1$,  рассматрим случайную   последовательность параметров $\bar s_j^{(q)}$, которая фигурирует в лемме 1.
 На этапе $j= q (mod\, k)$ положим  $\bar s_j=\bar s_j^{(q)}$, но
параметры $s_j(jp)$  модифицируются по правилу $(Sidon)$. Это  искомая конструкция
в случае, когда для допустимых функций коэффициенты удовлетворяют равенству $\sum_0^n c_k=1.$ 
Если $\sum_0^\infty  c_k=c  <1,$  параметры $s_j(jp)$ также модифицируются по правилу $(Sidon)$, но   при $cr_j<i\leq r_j$    все параметры $s_j(i)$ следует сделать сидоновскими. 
Будем называть полученные преобразования ранга один $(P_1, \dots, P_k)$-конструкциями. 

\vspace{2mm}
\bf Теорема 5. \it Пусть для   $(P_1, \dots, P_k)$-конструкции $T$ оператор $Q\neq 0$
лежит в слабом замыкании степеней $T^m$. Тогда найдутся последовательность ${m_p}$, для которой $T^{m_p}\to_w Q,$  натуральное $n$ и целые числа  $a_{1},a_{2},\dots, a_{n}, z$  такие, что для всякого $p$ для некоторых этапов $j_1>j_2>\dots>j_n$ выполнено  
$$m_p= a_{1}h_{j_1}+ a_{2}h_{j_2}+\dots+ a_{n}h_{j_n}+ z.$$
В этом случае оператор $Q$ имеет следующий вид:
$$Q=T^z \prod_{i=1}^n  Q_i^{|a_i|}, \ \  Q_i\in \{P_1(T),\dots, P_k(T), P_1(T^\ast),\dots, P_k(T^\ast)\}.$$\rm

\end{document}